# LES PROBABILITÉS DE DÉFAILLANCE COMME INDICATEURS DE PERFORMANCE DES BARRIÈRES TECHNIQUES DE SÉCURITÉ – APPROCHE ANALYTIQUE

Florent BRISSAUD & Brice LANTERNIER

*Institut National de l'Environnement Industriel et des Risques (INERIS)*
*Parc Technologique ALATA BP-2, 60550 Verneuil-en-Halatte, France*
*+33 3 44 55 69 89, florent.brissaud@ineris.fr & brice.lanternier@ineris.fr*

**Résumé :**

Le Code de l'environnement français impose aux industriels assujettis aux études de danger d'inclure des critères probabilistes, notamment pour déterminer le niveau de confiance des mesures de maîtrise des risques. Cet article présente les probabilités de défaillance en tant qu'indicateurs de performance des barrières techniques de sécurité. Une approche analytique est retenue. Des équations génériques sont alors proposées, permettant d'évaluer des probabilités de défaillance en intégrant différents paramètres : taux de défaillance, architecture de la barrière, tests de révision complets et partiels. Dans de nombreux cas, les résultats obtenus sont directement utilisables pour évaluer les niveaux de confiance.

**Abstract:**

French environmental laws require industrialists to include probability criteria in risk assessments, especially to define confidence levels for risk management measures. This paper presents the failure probabilities as efficient indicators for technical safety barrier performances. Generic formulas are proposed to evaluate these probabilities, including failure rate, barrier architecture, full and partial proof tests. In many cases, these results can be directly used to assess safety barrier confidence levels.

**Mots clés :** Barrière Technique de Sécurité, Système Instrumenté de Sécurité, SIS, Probabilité de défaillance, PFD, PFH, Niveau d'intégrité de sécurité, SIL, IEC 61508, Niveau de confiance

**Keywords:** Safety Technical Barrier, Safety Instrumented System, SIS, Probability of failure, PFD, PFH, Safety Integrity Level, SIL, IEC 61508, Confidence level

## 1. Introduction et contexte réglementaire

En France, le plan de prévention des risques technologiques (PPRT) repose principalement sur la réglementation des installations classées, régie par le Code de l'Environnement [1]-Livre V. Les installations classées pour la protection de l'environnement (ICPE), soumises à autorisation (autorisation (A) ou autorisation avec servitude (AS)), doivent fournir une étude de danger (EDD) [1]-Article L551. Elle a pour but de caractériser, analyser, évaluer, prévenir et réduire les risques de l'installation, ainsi que préciser l'ensemble des mesures de maîtrise des risques (MMR) mises en œuvre. Après la catastrophe d'AZF à Toulouse en Septembre 2001, la loi du 30 Juillet 2003 [2], dite « loi Bachelot », a modifié le Code l'Environnement en introduisant des critères probabilistes [1]-Article L512-1, alinéas 4 et 5 :

■ l'étude de danger donne lieu à *une analyse de risques qui prend en compte la probabilité d'occurrence, la cinétique et la gravité des accidents potentiels selon une méthodologie qu'elle explicite*
■ elle définit et justifie *les mesures propres à réduire la probabilité et les effets de ces accidents*

Les règles minimales relatives à l'évaluation et à la prise en compte de la probabilité d'occurrence sont déterminées par l'Arrêté du 29 Septembre 2005 [3]. Bien que la méthode soit libre, celle-ci doit être pertinente, mener à des résultats cohérents et répondre au principe de proportionnalité [4]. L'évaluation probabiliste du risque peut s'appuyer sur [3]-Article 2, alinéa 2 :

■ *la fréquence des évènements initiateurs spécifiques ou génériques*

■ *les niveaux de confiance des mesures de maîtrise des risques agissant en prévention ou en limitation des effets* (i.e. les barrières de sécurité)

Trois critères de performance d'une barrière de sécurité sont retenus dans la Circulaire du 7 Octobre 2005 [4] : efficacité, temps de réponse et niveau de confiance. Pour une efficacité et un temps de réponse donnés, c'est le niveau de confiance qui introduit le critère de probabilité : *le niveau de confiance est l'architecture et la classe de probabilité (…) pour qu'une barrière, dans son environnement d'utilisation, assure la fonction de sécurité pour laquelle elle a été choisie* [4]. Il s'agit ainsi d'une vulgarisation de la notion de disponibilité telle que définie par la CEI 60050 [5] : *la disponibilité est l'aptitude d'une entité à être en état d'accomplir une fonction requise dans des conditions données, à un instant donné ou* [en moyenne] *pendant un intervalle de temps donné*. À noter que le déclenchement intempestif d'une barrière de sécurité fait parti des critères non retenus par la réglementation, bien que cela puisse parfois conduire à des situations de danger.

Ici, nous traiterons exclusivement des barrières techniques de sécurité (BTS) assimilables à des systèmes instrumentés de sécurité (SIS). En Europe, la principale référence pour l'évaluation de ces barrières est la CEI 61508 [6]. Chaque fonction de sécurité allouée à la barrière est spécifiée en termes de niveau d'intégrité de sécurité (SIL) [6]-Partie 1, Section 7.6. Un SIL peut aller de 1 pour le niveau de sécurité le plus faible, à 4 pour le plus élevé [6]-Partie 4, Section 3.5.6. La Circulaire du 7 Octobre 2005 précise qu'il est possible de déterminer le niveau de confiance d'une barrière par l'intermédiaire du SIL [4] (la réciproque est fausse). Pour la suite, nous nous focaliserons donc sur le SIL, plus adapté aux SIS, et dont les exigences sont définies grâce à la CEI 61508 [6]. De nombreux critères qualitatifs, notamment sur l'architecture fonctionnelle de la barrière, sont alors requis pour un SIL donné [6]. La partie quantitative se fait quant à elle par l'un des deux indicateurs suivants [6]-Partie 1, Section 7.6.2.9 :

■ probabilité de défaillance à la demande (*PFD*) pour une barrière à faible demande de sollicitation
■ probabilité de défaillance par heure (*PFH*) pour une barrière à forte demande de sollicitation

Une barrière de sécurité est en mode de sollicitation à faible demande lorsque la fréquence des demandes d'opération n'est pas plus grande que une par an et pas plus grande que le double de la période des tests de révision [6]-Partie 4, Sections 3.5.12 et 13. Dans les autres cas, le mode de sollicitation est à forte demande. Le mot « défaillance » fera ici uniquement référence aux défaillances dangereuses et non détectées en temps réel. La correspondance entre *PFD* ou *PFH* et SIL est donnée dans le Tableau 1.

***Tableau 1.*** *SIL avec valeur cible de défaillance [6]-Partie 1, Section 7.6.2.9*

| SIL | Mode de sollicitation à faible demande | Mode de sollicitation à forte demande |
|---|---|---|
| 4 | $10^{-5} \leq PFD < 10^{-4}$ | $10^{-9} \leq PFH < 10^{-8}$ |
| 3 | $10^{-4} \leq PFD < 10^{-3}$ | $10^{-8} \leq PFH < 10^{-7}$ |
| 2 | $10^{-3} \leq PFD < 10^{-2}$ | $10^{-7} \leq PFH < 10^{-6}$ |
| 1 | $10^{-2} \leq PFD < 10^{-1}$ | $10^{-6} \leq PFH < 10^{-5}$ |

Dans cet article, nous discuterons tout d'abord de l'utilisation de critères probabilistes en tant qu'indicateurs de performance pour les barrières techniques de sécurité. Ensuite nous proposerons quelques équations permettant d'évaluer ces probabilités, en accord avec la CEI 61508 [6], et répondant aux exigences de la réglementation française sur les installations classées [1]. L'approche analytique possède en particulier l'avantage de fournir des résultats génériques pour certaines configurations données, ainsi que de pouvoir s'adapter facilement aux données d'entrée.

## 2. Les Probabilités comme indicateurs de performance des barrières de sécurité

### 2.1. Quelle interprétation prendre des probabilités ?

Il existe deux grandes catégories d'interprétations des probabilités :

■ les interprétations fréquentistes (dites également « probabilités physiques » ou « objectives »)
■ les interprétations Bayésiennes (dites également « probabilités de l'évidence » ou « subjectives »)

La première n'a de sens que lorsqu'une expérience aléatoire existe. La probabilité d'un évènement est alors définie par la limite de sa fréquence relative d'occurrence, lorsque l'expérience est répétée un très grand nombre de fois. Pour l'évaluation de la probabilité de défaillance d'une barrière de sécurité, cette approche est souvent implicitement retenue lorsque des bases de données sont utilisées. L'hypothèse est alors faite que la barrière étudiée est physiquement similaire et présente dans les mêmes conditions d'exploitation que celles recensées dans les bases. Les données collectées sont ainsi supposées être issues d'une même expérience aléatoire. Dans la pratique, de nombreux facteurs d'influence impactent la fiabilité de chaque barrière, rendant la population étudiée souvent très hétérogène [7] et faussant ainsi la notion d'expérience aléatoire. La définition fréquentiste des probabilités n'est alors plus appropriée.

L'interprétation Bayésienne [NDA : il est ici fait référence à une interprétation des probabilités et non aux inférences statistiques] peut quant à elle être utilisée pour n'importe quelle assertion, même lorsqu'aucune expérience aléatoire n'est admise. Une probabilité représente ici le degré de croyance (i.e. la mesure de l'état des connaissances) en la réalisation d'un évènement. Cet évènement peut être très rare, voir ne jamais avoir été observé auparavant, pour autant qu'il existe néanmoins des raisons de croire qu'il peut se produire. Considérons l'évènement « une défaillance de la barrière de sécurité se produira dans l'année ». Une fois l'année écoulée, il se serra alors possible de constater si oui ou non cette assertion aura été vérifiée. Ainsi, une « probabilité exacte » n'aurait pu être que *0* (non réalisation) ou *1* (réalisation). L'attribution d'une valeur strictement comprise entre *0* et *1*, bien qu'implicitement « fausse », représente initialement le degré de confiance en la réalisation de l'évènement, compte tenu des informations disponibles. L'utilisation d'une probabilité ne signifie donc pas qu'une vérité chiffrée a été établie mais qu'une croyance a été traduite sur une échelle allant de *0* à *1*. À noter que ce degré de croyance de l'interprétation Bayésienne peut s'appuyer par exemple sur des données issues d'approches fréquentistes.

### 2.2. Quel indicateur de performance retenir ?

Un parallèle est observable entre la définition Bayésienne d'une probabilité et l'utilisation des termes « niveau de confiance » par la législation française. L'évaluation d'une telle probabilité de défaillance semble ainsi être l'indicateur adéquat à la détermination d'un niveau de confiance. Bien que, par nature, la valeur obtenue ne puisse être « exacte » (sauf à prendre la valeur *0* ou *1*), nous pouvons néanmoins définir deux qualités importantes pour un tel indicateur :

■ cohérence (terme utilisé dans la Circulaire du 7 Octobre 2005 [4]) i.e. la probabilité de défaillance obtenue pour une barrière doit être supérieure à celle obtenue pour une seconde barrière si elle est jugée a priori « moins fiable » (i.e. les données d'entrée sont globalement jugées « moins favorables »)
■ robustesse i.e. dans la mesure du possible, la méconnaissance sur les données exploitées pour l'évaluation de la probabilité ne doit laisser qu'un minimum de doutes sur la cohérence des résultats

D'après la Circulaire du 7 Octobre 2005, en plus du niveau de confiance, deux autres critères sont considérés : l'efficacité et le temps de réponse [4]. L'efficacité est relative aux principes de dimensionnement adapté et de résistance aux contraintes spécifiques [4]. Une inefficacité de la barrière correspond ainsi à une défaillance systématique, c'est-à-dire *liée d'une manière certaine à une cause qui ne peut être éliminée que par une modification de la conception, du procédé de fabrication (…) ou d'autres facteurs appropriés* [5, 6]. L'efficacité d'une barrière peut ainsi être intégrée à la probabilité de défaillance grâce à la prise en compte des probabilités de défaillance systématique (*PSF*) [6, 8].

Le temps de réponse doit quant à lui être spécifié dans la définition de la fonction à évaluer. Une fonction est implicitement assurée si elle est pleinement effectuée dans les temps qui lui sont impartis. De plus, il est parfois judicieux de prendre en compte le délai aléatoire d'exécution de la fonction de sécurité d'une barrière dans l'évaluation de la probabilité de défaillance.

Une probabilité de défaillance permet ainsi de regrouper dans le même indicateur les critères réglementaires de performance d'une barrière de sécurité. La nature d'un tel indicateur est alors avantageuse pour l'évaluation globale des risques et de la gestion des mesures de maîtrise des risques.

### 2.3. *Le choix d'un modèle d'évaluation des probabilités de défaillance*

Différents modèles existent pour évaluer une probabilité de défaillance (*PFD* ou *PFH*) : blocs diagrammes de fiabilité, arbres de défaillance, graphes de Markov, Réseaux de Petri… Ici, nous utiliserons directement une approche analytique, car relativement souple en termes d'hypothèses et de paramètres exploités. De plus, les équations de fiabilité offrent des résultats génériques pour un ensemble de barrières de sécurité présentant des caractéristiques semblables. L'analyse des résultats obtenus en fonction des paramètres d'entrée, par exemple à des fins d'optimisation, est également plus directe. Voici certains paramètres qu'il est envisageable d'inclure :

■ les propriétés générales de la barrière, représentées par les taux de défaillance. Cette notion purement mathématique est elle-même le résultat d'un modèle. Le retour d'expérience est souvent utilisé pour le choix des taux de défaillance, participant ainsi de façon significative à l'évaluation du SIL. De nombreux facteurs d'influence peuvent également être inclus par l'intermédiaire de ce paramètre [7]

■ l'architecture (terme utilisé dans la Circulaire du 7 Octobre 2005 [4]) i.e. les propriétés de redondance des éléments de la barrière, souvent symbolisées par les architectures *MooN* (cf. Section 3.1)

■ les instants, périodes et durées des tests de révision (spécifiés dans la CEI 61508 [6]). Ces tests permettent de détecter les défaillances non automatiquement détectables en temps réel. Ils peuvent être considérés comme complets ou partiels (i.e. détection d'une partie uniquement des défaillances)

■ la durée des réparations (allant du diagnostic à la remise en état de fonctionnement) de la barrière lorsqu'une défaillance a été détectée. Si des mesures adéquates sont prises pour maintenir la fonction de sécurité avec le même SIL, ou bien que la situation de danger est écartée (e.g. arrêt du processus) pendant la durée de réparation, alors ce paramètre peu rationnellement être omis du calcul de *PFD* ou *PFH*

■ les causes communes de défaillance, qui atténuent généralement l'intérêt des redondances matérielles au sein d'une barrière de sécurité. En fonction du modèle utilisé, différents paramètres peuvent être exploités comme par exemple le *facteur β* ou le $C_{MooN}$ [8]

De nombreux autres paramètres peuvent être pris en compte dans l'évaluation des probabilités de défaillance. Leurs choix doivent être motivés par la quantité et la qualité des informations disponibles. Par exemple, il est possible d'exploiter un taux de défaillance à plusieurs paramètres (e.g. variable dans le temps). La question est alors de savoir si le nombre de données est suffisant pour prendre en compte cette particularité dans l'évaluation, ou bien si l'hypothèse d'un taux de défaillance constant est raisonnable.

## 3. Approche analytique pour l'expression des probabilités de défaillance

### 3.1. *Notations et hypothèses*

■ *MooN* – **architecture de la barrière** : la barrière est constituée de *N* éléments et est opérante si et seulement si au moins *M* de ces éléments sont opérants (avec $M \leq N$)

■ *λ* – **taux de défaillance [dangereuse et non détectée en temps réel] de chacun des *N* éléments de la barrière** [par heure] : sachant que l'élément est opérant à l'instant *t*, $\lambda \cdot \Delta t$ est la probabilité pour qu'une défaillance de l'élément se produise dans l'intervalle *[t, t+Δt]* lorsque *Δt* tend vers *0*.

■ *$\lambda_b(t)$* – **taux de défaillance [dangereuse et non détectée en temps réel] de la barrière** [par heure] : sachant que la barrière est opérante à l'instant *t*, $\lambda_b(t) \cdot \Delta t$ est la probabilité pour qu'une défaillance de la barrière se produise dans l'intervalle *[t, t+Δt]* lorsque *Δt* tend vers *0*

■ *R(t)* – **fonction de fiabilité de la barrière** : probabilité pour qu'aucune défaillance de la barrière ne se produise dans l'intervalle *[0, t]*, l'instant *0* correspondant à la dernière remise à neuf de la barrière

■ *$T_1$* – **période des tests de révision complets** [heure] : après chaque test de révision complet, la barrière de sécurité est comme neuve (les tests sont parfaits). Lorsqu'aucun test de révision complet n'est envisageable, la remise à neuf peut consister en un remplacement total de la barrière. La durée des tests (e.g. tests en ligne) et des actions de maintenance ne sera pas prise en compte (cf. Section 2.3)

■ *$T_0$* – **période des tests de révision partiels** [heure] : après chaque test de révision partiel, seule une proportion *E* des défaillances de chacun des éléments de la barrière est détectée. Entre chaque test de révision complet, *n* tests de révision partiels sont effectués à intervalle régulier (i.e. $T_0 = T_1/n$)

■ *PFD(t)* – **probabilité de défaillance à la demande de la barrière à l'instant *t*** [i.e. indisponibilité instantanée] : sous les hypothèses précédentes, on vérifie *PFD(t) = 1-R(t)* dans l'intervalle *[0, T₁]*
■ *PFD* – **probabilité moyenne de défaillance à la demande de la barrière** [i.e. indisponibilité moyenne], [NDA : souvent notée *PFD_{avg}*] : d'après les hypothèses, calculée sur l'intervalle *[0, T₁]*
■ *PFH(t)* – **probabilité de défaillance par heure de la barrière à l'instant *t***
■ *PFH* – **probabilité moyenne de défaillance par heure de la barrière**

### 3.2. Probabilité de défaillance à la demande (PFD), cas basique sans test de révision partiel

D'après les notations et les hypothèses de la Section 3.1., la fonction de fiabilité de la barrière est [9] :

$$R(t) = \sum_{k=M}^{N} \binom{N}{k} \cdot e^{-k \cdot \lambda \cdot t} \cdot \left(1 - e^{-\lambda \cdot t}\right)^{N-k} \quad \text{pour } t \in [0, T_1] \quad (1)$$

On en déduit la probabilité de défaillance à la demande à l'instant *t*, et moyenne [10] :

$$PFD(t) = 1 - \sum_{x=M}^{N} S(M,N,x) \cdot e^{-x \cdot \lambda \cdot t} \quad \text{pour } t \in [0, T_1] \quad (2)$$

$$PFD = 1 - \sum_{x=M}^{N} S(M,N,x) \cdot \frac{1 - e^{-x \cdot \lambda \cdot T_1}}{x \cdot \lambda \cdot T_1} \quad (3)$$

Avec la somme suivante, indépendante du temps *t* :

$$S(M,N,x) = \sum_{k=M}^{x} \binom{N}{x} \cdot \binom{x}{k} \cdot (-1)^{x-k} \quad \text{pour } x = M,...,N \quad (4)$$

La Figure 1 représente *PFD(t)* et *PFD* pour une configuration donnée. Les valeurs de la somme *S(M,N,x)* sont données dans le Tableau 2 pour différentes architectures *MooN*. Par exemple, pour une barrière d'architecture *2oo3*, les valeurs correspondantes sont *S(2,3,2) = 3* et *S(2,3,3) = -2*, c'est-à-dire que :

$$PFD(t) = 1 - 3 \cdot e^{-2 \cdot \lambda \cdot t} + 2 \cdot e^{-3 \cdot \lambda \cdot t} \quad \text{pour } t \in [0, T_1] \; ; \; PFD = 1 - 3 \cdot \frac{1 - e^{-2 \cdot \lambda \cdot T_1}}{2 \cdot \lambda \cdot T_1} + 2 \cdot \frac{1 - e^{-3 \cdot \lambda \cdot T_1}}{3 \cdot \lambda \cdot T_1} \quad (5)\text{-}(6)$$

***Tableau 2.*** *S(M,N,x) sous la forme S(M,N,M) / S(M,N,M+1) / ... / S(M,N,N)*

|  | *N = 1* | *N = 2* | *N = 3* | *N = 4* |
|---|---|---|---|---|
| *M = 1* | 1 | 2 / -1 | 3 / -3 / 1 | 4 / -6 / 4 / -1 |
| *M = 2* |  | 1 | 3 / -2 | 6 / -8 / 3 |
| *M = 3* |  |  | 1 | 4 / -3 |
| *M = 4* |  |  |  | 1 |

Lorsque *λ·t* est petit (*λ·t << 10⁻²*), le développement limité d'ordre *1* de la fonction exponentielle au voisinage de *0* nous donne l'approximation conservative *1-e^{-λ·t} ≈ λ·t*. Cette démarche nous permet d'obtenir les équations approchées suivantes, lorsque *λ·T₁* est petit (*λ·T₁ << 10⁻²*) [10] :

$$PFD(t) \approx \binom{N}{M-1} \cdot (\lambda \cdot t)^{N-M+1} \quad \text{pour } t \in [0, T_1] \; ; \; PFD \approx \binom{N}{M-1} \cdot \frac{(\lambda \cdot T_1)^{N-M+1}}{N-M+2} \quad (7)\text{-}(8)$$

Les équations (7) et (8) peuvent être obtenues par un raisonnement sur les plus courts chemins de coupe [9]-p. 432. Elles sont d'autant plus proches des équations (2) et (3) que la valeur de *T₁* est petite. Par exemple, dans le cas d'un système à un seul élément (i.e. architecture *1oo1*), on obtient *PFD ≈ λ·T₁/2*.

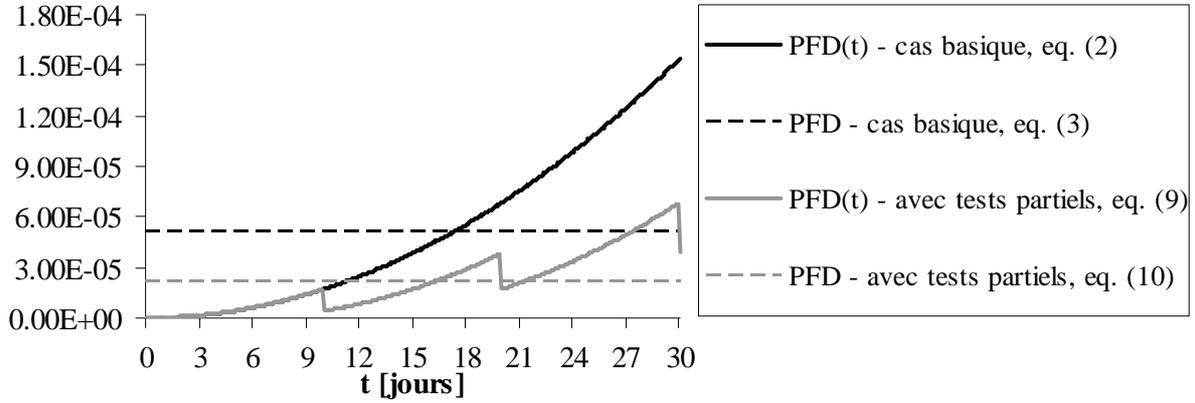

**Figure 1.** *PFD(t) et PFD pour le cas basique sans test de révision partiel (équations (2) et (3))
et pour le cas avec tests de révision partiels (équations (9) et (10))
architecture 2oo3, $\lambda = 1 \cdot 10^{-5}$ heure$^{-1}$, $T_1$ = 30 jours, E = 0.5, $T_0$ = 10 jours, $n = T_1/T_0 = 3$*

### 3.3. Probabilité de défaillance à la demande (PFD) avec tests de révision partiels

Les tests de révision partiels (i.e. incomplets) ne permettent de détecter qu'une certaine proportion des défaillances. Plusieurs raisons peuvent justifier le recours à ce type de tests :

■ les tests complets doivent généralement être « physiques » (e.g. excitation des organes sensibles d'un capteur), ce qui est couteux. Ils sont alors souvent remplacés par des tests purement « électroniques » (e.g. simulation électronique d'une sollicitation) qui ne permettent pas de couvrir toutes les défaillances
■ les tests complets impliquent souvent des arrêts de production (e.g. coupure de l'alimentation électrique par un arrêt d'urgence, coupure d'un flux par une vanne de sécurité), souvent jugés inacceptables. Des tests partiels (e.g. fermeture d'un quart de tours d'une vanne) sont alors préférés
■ pour certaines barrières de sécurité, un test complet ne peut s'envisager sans la dégradation ou la destruction de celle-ci (e.g. systèmes à rupture)
■ seul un test en situation réel pourrait prétendre à être complet. Or, cela provoquerait dans de nombreux cas plus de dangers que de prévention (e.g. détection d'incendie ou de gaz toxique, seuil de surpression)

D'après les notations et les hypothèses de la Section 3.1, la probabilité de défaillance à la demande à l'instant *t*, et moyenne, sont dans ce cas [10] :

$$PFD(t) = 1 - \sum_{x=M}^{N} S(M,N,x) \cdot e^{x \cdot E \cdot \lambda \cdot \left\lfloor \frac{t}{T_0} \right\rfloor \cdot T_0} \cdot e^{-x \cdot \lambda \cdot t} \quad \text{pour } t \in [0, T_1] \tag{9}$$

$$PFD = 1 - \sum_{x=M}^{N} S(M,N,x) \cdot T(n, E, \lambda, T_0, x) \cdot \frac{1 - e^{-x \cdot \lambda \cdot T_0}}{x \cdot \lambda \cdot T_0} \tag{10}$$

Avec $\lfloor \ \rfloor$ le symbole de la partie entière, $S(M,N,x)$ définie par l'équation (4), et la somme suivante :

$$T(n, E, \lambda, T_0, x) = \frac{1}{n} \cdot \sum_{p=0}^{n-1} e^{-x \cdot (1-E) \cdot \lambda \cdot p \cdot T_0} \tag{11}$$

La Figure 1 représente *PFD(t)* et *PFD* suivant les équations (9) et (10), pour une configuration donnée. Lorsque $\lambda \cdot T_1$ est petit ($\lambda \cdot T_1 << 10^{-2}$), nous pouvons obtenir les équations approchées suivantes [10] :

$$PFD(t) \approx \binom{N}{M-1} \cdot \left( \lambda \cdot t - E \cdot \lambda \cdot \left\lfloor \frac{t}{T_0} \right\rfloor \cdot T_0 \right)^{N-M+1} \quad \text{pour } t \in [0, T_1] \tag{12}$$

$$PFD \approx \binom{N}{M-1} \cdot \frac{(\lambda \cdot T_0)^{N-M+1}}{N-M+2} \cdot V(M,N,n,E) \tag{13}$$

Avec la somme suivante :

$$V(M,N,n,E) = \frac{1}{n} \cdot \sum_{p=0}^{n-1} \left[ (1 + p \cdot (1-E))^{N-M+2} - (p \cdot (1-E))^{N-M+2} \right] \tag{14}$$

### 3.4. Probabilité de défaillance par heure (PFH)

Bien que la *PFH* soit établie comme étant une probabilité [6]-Parties 1 et 4, la CEI 61508 l'assimile tantôt à une « fréquence des défaillances » (i.e. intensité inconditionnelle des défaillances), tantôt à un « taux de défaillance » [6]-Partie 4, Section 7.6.2.9, Note 4. Ces deux notions étant des fréquences (dimension [temps]$^{-1}$ et éventuellement plus grandes que *1*), et non des probabilités (sans dimension et comprise entre *0* et *1*), elles introduisent certaines confusions. De plus, la fréquence des défaillances tend vers *0* pour des systèmes non réparés (i.e. une fois le système défaillant, une remise en état est nécessaire avant une prochaine défaillance), ce qui est le cas dans l'intervalle *[0 ; T₁]* lorsque les défaillances ne sont pas détectées. Cela peut donc conduire à une fréquence moyenne des défaillances qui est décroissante en fonction de *T₁*. Un taux de défaillance peut également être décroissant par périodes (e.g. pour un système *1oo2* avec des taux de défaillance différents [9]-p.56, pour un système en période de rodage etc.). Dans certains cas, la *PFH* alors obtenue par l'une ou l'autre de ces deux approches peut alors être réduite, et donc amener à un meilleur SIL, en augmentant les périodes de tests de révision, ce qui est un non-sens.

Pour lever cet ambiguïté, nous proposons d'interpréter la notion de *PFH* comme la *probabilité de défaillance en une heure* : sachant que la barrière de sécurité est opérante à l'instant *t*, *PFH(t)* est la probabilité pour qu'une défaillance de la barrière se produise dans l'intervalle *[t, t+1 heure]*. Il est ainsi fait référence à la fonction de défiabilité conditionnelle à l'instant *t*, pour une durée d'une heure [10] :

$$PFH(t) = 1 - \frac{R(t+1\ heure)}{R(t)} = 1 - \frac{1 - PFD(t+1\ heure)}{1 - PFD(t)} \quad \text{pour } t \in [0, T_1] \tag{16}$$

$$PFH = 1 - \frac{1}{T_1} \cdot \int_0^{T_1} \left( \frac{R(t+1\ heure)}{R(t)} \right) \cdot dt = 1 - \frac{1}{T_1} \cdot \int_0^{T_1} \left( \frac{1 - PFD(t+1\ heure)}{1 - PFD(t)} \right) \cdot dt \tag{17}$$

Définie ainsi, la *PFH* est une probabilité (sans dimension, résultat compris entre *0* et *1*). Si le taux de défaillance de la barrière présente peu de variations dans le temps par rapport à *1 heure* (i.e. $|\lambda_b(t+1\ heure) - \lambda_b(t)|$ est faible pour tout *t*), alors la probabilité de défaillance en une heure à l'instant *t* peut être approchée par le taux de défaillance de la barrière, multiplié par *1 heure* (le résultat reste alors sans dimension). Cette particularité explique que les équations approchées (18) et (19) sont semblables à celles que l'on peut déduire d'une définition par fréquence ou par taux de défaillance [10] :

$$PFH(t) \approx \lambda_b(t) \cdot 1\ heure \quad \text{avec } \lambda_b(t) \approx \binom{N}{M} \cdot M \cdot \lambda^{N-M+1} \cdot t^{N-M} \quad \text{pour } t \in [0, T_1] \tag{18}$$

$$PFH \approx \binom{N}{M-1} \cdot \lambda^{N-M+1} \cdot T_1^{N-M} \cdot 1\ heure \tag{19}$$

On remarque alors la relation approchée suivante, entre la *PFH* et la *PFD* :

$$PFH \approx \frac{PFD(T_1)}{T_1} \cdot 1\ heure \approx PFD \cdot \frac{N-M+2}{T_1} \cdot 1\ heure \tag{20}$$

## 4. Conclusion

Nous avons présenté dans cet article comment la définition d'une probabilité de défaillance, notamment d'après l'interprétation Bayésienne que l'on peut en faire, permet d'intégrer dans un même indicateur les critères réglementaires de performance d'une barrière technique de sécurité. En accord avec la réglementation française sur les études de danger, ainsi que la CEI 61508, des équations génériques pour le calcul des probabilités de défaillance à la demande (*PFD*), ou par heure (*PFH*), ont été proposées. Plusieurs paramètres peuvent alors être pris en compte : taux de défaillance, architecture de la barrière, tests de révision complets et partiels. Une formulation de la *PFH* a également été proposée, cherchant à traduire cet indicateur en probabilité et non en fréquence. Dans de nombreux cas, les résultats obtenus sont directement exploitables pour évaluer le niveau de confiance d'une barrière technique de sécurité.